\documentclass{amsart}

\def\1#1{\overline{#1}}
\def\2#1{\widetilde{#1}}
\def\3#1{\widehat{#1}}
\def\4#1{\mathbb{#1}}
\def\5#1{\mathfrak{#1}}
\def\6#1{{\mathcal{#1}}}

\def\Z{{\4Z}}

\def\Label#1{\label{#1}}

\def\aut{{\sf aut}}
\def\Aut{{\sf Aut}}

\def\Re{{\sf Re}\,}
\def\Im{{\sf Im}\,}

\begin{document}

\def\codim{{\rm codim}}
\def\crd{\dim_{{\rm CR}}}
\def\crc{{\rm codim_{CR}}}

\def\phi{\varphi}
\def\eps{\varepsilon}
\def\d{\partial}
\def\a{\alpha}
\def\b{\beta}
\def\g{\gamma}
\def\G{\Gamma}
\def\D{\Delta}
\def\Om{\Omega}
\def\k{\kappa}
\def\l{\lambda}
\def\L{\Lambda}
\def\z{{\bar z}}
\def\w{{\bar w}}
\def\Z{{\mathbb Z}}
\def\t{\tau}
\def\th{\theta}

\emergencystretch15pt
\frenchspacing

\newtheorem{Thm}{Theorem}[section]
\newtheorem{Cor}[Thm]{Corollary}
\newtheorem{Pro}[Thm]{Proposition}
\newtheorem{Lem}[Thm]{Lemma}

\theoremstyle{definition}\newtheorem{Def}[Thm]{Definition}

\theoremstyle{remark}
\newtheorem{Rem}[Thm]{Remark}
\newtheorem{Exa}[Thm]{Example}
\newtheorem{Exs}[Thm]{Examples}

\def\bl{\begin{Lem}}
\def\el{\end{Lem}}
\def\bp{\begin{Pro}}
\def\ep{\end{Pro}}
\def\bt{\begin{Thm}}
\def\et{\end{Thm}}
\def\bc{\begin{Cor}}
\def\ec{\end{Cor}}
\def\bd{\begin{Def}}
\def\ed{\end{Def}}
\def\br{\begin{Rem}}
\def\er{\end{Rem}}
\def\be{\begin{Exa}}
\def\ee{\end{Exa}}
\def\bpf{\begin{proof}}
\def\epf{\end{proof}}
\def\ben{\begin{enumerate}} 
\def\een{\end{enumerate}}
\def\beq{\begin{equation}}
\def\eeq{\end{equation}}

\newcommand{ \ssf}{F^*}
\newcommand{ \jj}{\iota}
\newcommand{\alb}{\bar a_l}
\newcommand{\pz}{P(z,\bar z)}
\newcommand{\bet}{\beta}
\newcommand{\ga}{\gamma}
\newcommand{\de}{\delta}
\newcommand{\ka}{\kappa}
\newcommand{\ab}[1]{\vert z\vert^{#1}}
\newcommand{\cdva}{{\mathbb C^2}}
\newcommand{\om}{\Omega}
\newcommand{\ml}{M^{k,l}_a}
\newcommand{\te}{\theta}
\newcommand{\kl}{\dfrac k{l^2-k}}
\newcommand{\klv}{\dfrac k{\vert l^2-k\vert}}
\newcommand{\kldlo}{\dfrac{(4k-l^2-4)k^2}{(4k-4)(k^2-l^2)}}
\newcommand{ \ffn}{F^{\nu}}
\newcommand{ \fn}{f^{\nu}}
\newcommand{\pl}{P^{k,l}_a}
\newcommand{ \ksi}{\xi}
\newcommand{\zz}{(z,\bar z)}
\newcommand{ \nl}{\newline}
\newcommand{\fz}{{F(z,\bar z,u)}}
\newcommand{\fzs}{{\textstyle{F^*(z^*,\bar z^\vert (\alpha,\hat \al) \vert_{\La}*,u^*)}}}
\newcommand{\dom}{{\textstyle{b\Omega \ }}}
\newcommand{ \al}{\alpha}
\newcommand{\ppm}{m'}
\newcommand{\db}{\overline{\partial}}
\newcommand{\tf}{\tilde F}
\newcommand{\la}{\lambda}
\newcommand{\La}{\Lambda}
\newcommand{\kkk}{\frac{k}2 \frac{k}2}

\newtheorem{theorem}{Theorem}[section]
\newtheorem{lemma}[theorem]{Lemma}

\theoremstyle{definition}
\newtheorem{definition}[theorem]{Definition}
\newtheorem{example}[theorem]{Example}
\newtheorem{xca}[theorem]{Exercise}

\theoremstyle{remark}
\newtheorem{remark}[theorem]{Remark}

\numberwithin{equation}{section}

\def\to{\rightarrow}
\def\aut{\operatorname{aut}}
\def\Aut{\operatorname{Aut}}
\def\Hom{\operatorname{Hom}}
\def\End{\operatorname{End}}
\def\d{\partial}
\def\Re{{\sf Re}\,}
\def\Im{{\sf Im}\,}

\title{Infinitesimal symmetries of weakly pseudoconvex manifolds}

\author{Shin-Young Kim}
\address{Institut Fourier, Universit\'e Grenoble Alpes}
\curraddr{}
\email{Shinyoung.Kim@univ-grenoble-alpes.fr}
\thanks{}
\author{Martin Kol\'a\v r}
\address{Department of Mathematics and Statistics, Masaryk University}
\curraddr{}

\email{mkolar@math.muni.cz}
\thanks{The authors were supported by the GACR grant GA17-19437S.}

\subjclass[2010]{Primary }

\keywords{}

\date{}

\dedicatory{}

\begin{abstract}
We classify the Lie algebras of infinitesimal CR automorphisms of weakly pseudoconvex hypersurfaces of finite multitype in $\mathbb C^N$. In particular,
we prove that such manifolds admit neither nonlinear rigid automorphisms, nor real or nilpotent rotations. As a consequence, this leads to a proof of a sharp 2-jet determination 
result for local automorphisms. 
Moreover, for  hypersurfaces which are not balanced, CR automorphisms are uniquely determined by their 1-jets.
The same classification is derived also for special models, given by sums of squares of  polynomials. 
In particular, in the case of homogeneous polynomials the Lie algebra of infinitesimal CR automorphisms is always three  graded.
The results provide an important necessary step for solving the local equivalence problem on weakly pseudoconvex manifolds.
\end{abstract}

\maketitle

\section{Introduction}

The study of weakly pseudoconvex manifolds goes back to the 
pioneering work of J. J. Kohn \cite{K}, in which he defined type of a point $p \in M \subseteq \mathbb C^2$ as the  lowest order (integer valued)
local CR invariant.

More refined rational local invariants of pseudoconvex boundaries in  $\mathbb C^N$  turned out crucial for characterizing subellipticity of the $\bar{\partial}$ - Neumann problem, 
as follows from the work of D. Catlin  (\cite{C, C1}).

Recent work of D. Zaitsev and Sung-Yeon Kim  shows that qualitative understanding of higher order CR invariants is also essential to understand effective termination of the Kohn algorithm. 
In relation to the work \cite{KZ, Z},  Dmitri Zaitsev
posed a question on possible symmetries of weakly psedoconvex manifolds of finite Catlin multitype. Can the
structure results of \cite{KM14} be improved under the pseudoconvexity assumption?

The main result of this paper shows that the structure of $\5g=\aut(M_{P},0)$ for weakly pseudoconvex hypersurfaces is indeed substantially simpler compared to the general case. In particular, 
there exist no nonlinear rigid vector fields on such manifolds. 
As one consequence, this leads to a sharp 2 - jet determination result.
Moreover, 2-jets are needed precisely in the case of a balanced model. In the non balanced case we obtain 
determination 
of local automorphisms by their 1-jets.

We also consider the sum of squares models. In the case of homogeneous Levi degenerate sums of squares models the structure of $\aut(M_{P},0)$ is particularly simple, 
having precisely three nonvanishing components, 
where $\5g _{-1}$ and  $\5g _{1}$ are one dimensional and  $\5g _{0}$ is generated by the Euler field and a subalgebra of $\mathfrak u(n)$.

The ultimate motivation for our work comes from the Poincar\'e problem on local equivalence of real hypersurfaces in $\mathbb C^N$. Let us briefly recall some of the history 
and some recent developments, which motivated our work. In the classical case 
of Levi nondegenerate manifolds, the problem was solved in the works of Cartan, Tanaka, Chern and Moser (\cite{ C32, C2, {CM74}, Po, T}. While Cartan, Tanaka and Chern applied 
differential geometric techniques, Moser developed a normal form approach, inspired by the normal form solution to the equivalence problem for analytic vector fields, originating 
also in the work of Poincar\'e. 

The case of singular Levi form presents completely new challenges, which are often more of algebraic than of differential-geometric nature.
A construction of a normal form for finite type Levi-degenerate hypersurfaces
in $\mathbb C^2$ , and the description of their symmetries, was given  by the second author in \cite{K06}.
In  combination with a convergence result of
Baouendi-Ebenfelt-Rothschild \cite{BER00}, this normal form 
solves the biholomorphic equivalence problem for this class. The normal form has been shown
to be convergent under some additional geometric conditions  by Kossovskiy and Zaitsev, \cite{KZ15}
but is divergent in general, as shown by Kol\'a\v r \cite{K12}.

In complex dimensions higher than two, local geometry of Levi degenerate
hypersurfaces is far  more complicated, even on the initial level.
For pseudoconvex hypersurfaces, D. Catlin (\cite{C}) introduced
a notion of multitype.
The entries of the Catlin multitype take rational values, but 
need not be integers,
anymore.
This approach provides a
defining equation of the form 
\begin{equation}
Im \;w = P(z, \bar z) + o_w(1), \label{pzz}
\end{equation}
where $P$ is a weighted homogeneous polynomial in the complex tangential variables 
$z = (z_1, \dots, z_n)$, 
with respect to the multitype weights $(\mu_1, \dots, \mu_n)$, $w$ is the normal variable and $ o_w(1)$ denotes terms of weight bigger than one.  Moreover, the multitype weights are the lexicographically smallest weights for which such a description is possible, hence providing a fundamental CR invariant.

We will denote by $M_P$ the corresponding polynomial model,
\begin{equation}\Label{model2}
M_{P}:=\{\Im w = P(z,\bar z)\}.
\end{equation}

The paper \cite{K10}  introduced an alternative 
approach to the Catlin multitype,
which  
can be applied to a general smooth hypersurface in $\mathbb
C^{n+1}$,  not necessarily pseudoconvex. Using this approach, it
proved biholomorphic equivalence of models,   and gave an explicit
description of biholomorphisms between different models. It also provided a constructive finite
algorithm for computing the multitype. 

Using this approach,  Kol\'a\v r, Meylan, Zaitsev showed in [KMZ14] that hypersurfaces of finite Catlin multitype provide the natural class of manifolds
for which a generalization of the Chern-Moser operator is well defined.

Using this operator, Kol\'a\v r, Meylan, Zaitsev proved that   the   Lie algebra of infinitesimal automorphisms
$\5g=\aut(M_{P},0)$ of $M_{P}$ admits the  weighted grading given
by
\begin{equation*}
\5g = \5g_{-1} \oplus \bigoplus_{j=1}^{n}\5g_{- \mu_j} \oplus \5g_{0}
\oplus \bigoplus_{\eta \in E}\5g_{\eta}
 \oplus \5g_{1},
 \label{muth}
\end{equation*}
\\[2mm]
where
$E$ is the set of integer combinations of the multitype weights, which lie between zero and one. 
They also obtained an explicit description of the graded components. As a consequence, they
proved that the automorphisms of M at p are uniquely determined by their {\it weighted 2-jets} at p.

Since the kernel of the generalized Chern-Moser operator corresponds to the Lie algebra $\aut (M_P, 0)$ of infinitesimal CR automorphisms of the model, this result 
gives a necessary tool for addressing the equivalence problem. However, full classification of such Lie algebras seems 
still unattainable. One of the main difficulties is the presence of the component (denoted below as  $\5g_c$), containg nonlinear rigid vector fileds, 
with arbitrarily high degree coefficients (only the weighted degree is controlled).

In this paper we show that in the most interesting case of pseudoconvex  Levi degenerate manifolds, the structure of  $\aut (M_P, 0)$,
is in fact much simpler, thus opening the possibility for a complete solution of the Poincar\'e equivalence problem in this class.

Throughout the paper, $M$ will denote a hypersurface of finite Catlin multitype, described by (1.1). 

We now formulate the main results of the paper. 

\bt  Assume that $M$ is pseudoconvex in a neighbourhood of $p$
and the associated model hypersurface $M_P$, given by \eqref{model2}
is holomorphically nondegenerate.

Then the   Lie algebra of infinitesimal automorphisms
$\5g=\aut(M_{P},0)$ of $M_{P}$ admits the  weighted grading given
by
\begin{equation}
\5g = \5g_{-1} \oplus \bigoplus_{j=1}^{n}\5g_{- \mu_j} \oplus \5g_{0} \oplus \5g_{\frac12} 
\oplus \5g_{1}.
\label{muth}
\end{equation}

\et

As a consequence, we obtain the following sharp jet determination result.

\bt 
 Assume that $M$ is pseudoconvex in a neighbourhood of $p$
and the associated model hypersurface \eqref{model2}
is holomorphically nondegenerate.
Then the
automorphisms of $M$  at $p$ are uniquely determined by their jets
of order $2$.
 Moreover, if $M_P$ is not balanced, the automorphisms of $M$  at $p$ are uniquely determined by their jets
of order $1$.
\et

Special domains have been used as a manageable test case in the context of the d-bar problem and effectivity of the Kohn algorithm (see e.g. 
\cite{S}, \cite{C}). 

The model coresponding to a special domain is a sum of squares models. Of course, it may happen that 
the model for a holomorphically nondegenerate hypersurface is holomorphically degenerate, a thus does not 
provide useful information concerning symmetries of the original manifold itself. 

For special models, given by sums of squares of homogeneous polynomials, we prove the following results. 
By $M_S$, we will denote a model given by
\begin{eqnarray}
\Im (w) = \sum_{j=1}^{k} |P_j(z)|^2.
\end{eqnarray}

\bt
Let  $M_{S}$ be the sum of squares-type homogeneous polynomial model of degree $ k > 2$. 
Then the   Lie algebra of infinitesimal automorphisms
$\5g=\aut(M_{S},0)$ of $M_{S}$ admits the weighted grading given
by
\begin{equation}
\5g = \5g_{-1} \oplus \5g_{ 0} \oplus \5g_{1}
\end{equation}
where $\5g_{-1}$ and $\5g_{1}$ are of real dimension one and $\5g_{0}$ is generated by the Euler field and a subalgebra of $\mathfrak u(n)$. 
\et

We obtain an analogous result in the case of a weighted homogeneous sum of squares model. Let us denote by $\kappa_M$ the number of multitype weights with 
$\mu_j = \frac12$.

\bt
Let  $M_{S}$ be the sum of squares-type weighted homogeneous polynomial model  of degree $ k > 2$. 
Then the   Lie algebra of infinitesimal automorphisms
$\5g=\aut(M_{S},0)$ of $M_{S}$ admits the weighted grading given
by
\begin{equation}
\5g = \5g_{-1} \oplus \5g_{-\frac12} \oplus\5g_{ 0} \oplus\5g_{\frac12 } \oplus \5g_{1},
\end{equation}
where $\5g_{-1}$ and $\5g_{1}$ are of real dimension one,  $\5g_{-\frac12}$ and  $\5g_{-\frac12} $ are of real dimension $2 \kappa_M$,  and $\5g_{0}$ is generated by the Euler field and a subalgebra of $\mathfrak u(n)$. 

\et

The paper is organized as follows. In Section 2 we recall the needed definitions and notation, which is used in the sequel. In Section 3 we consider rotations and 
prove that psedoconvex models do not admit any real or nilpotent rotations. Section 4 considers generalized rotation. We show that such symmetries cannot occur on weakly psedoconvex manifolds. 
In Section 5 we consider sum of squares models. We prove Theorem 1.3 and 1.4. 
Theorems 1.1 and 1.2 are proved in Section 6.

\section{Preliminaries}

In this section we recall the definitions and notation needed in the rest of the paper (for more details, see e.g. \cite{KM14}).

\begin{definition}Let $n \in \mathbb N$ be an integer. A \emph{weight} is an $n$-tuple of nonnegative rational numbers $\mu=(\mu_1, \cdots, \mu_n)$, 
where $0\leq \mu_j \leq \frac{1}{2}$ and $\mu_{j} \geq \mu_{j+1}$ such that there exist  an $n$-tuple of non-negative integers $(\alpha_1, \cdots ,\alpha_n)$ satisfying $\alpha_j \neq 0$ if $\mu_j \neq 0$ for each $j$ and\begin{eqnarray} \sum_{j}\alpha_j\mu_j=1.
\end{eqnarray}

If $\alpha=(\alpha_1, \cdots ,\alpha_n)$ is a multiindex, for the given weight $\mu=(\mu_1, \cdots, \mu_n)$, \emph{the weighted length of $\alpha$} is
\begin{eqnarray*}
|\alpha|_{\mu}=
\sum_{j=1}^{n}\alpha_j \mu_j.
\end{eqnarray*}
Similarly, if $\alpha=(\alpha_1, \cdots ,\alpha_n)$ and $\hat \alpha=(\hat \alpha_1, \cdots ,\hat \alpha_n)$ are two multiindices, 
for the given weight $\mu=(\mu_1, \cdots, \mu_n)$, \emph{the weighted length of $(\alpha,\hat \alpha)$} is
\begin{eqnarray*}
|(\alpha,\hat \alpha)|_{\mu}=\sum_{j=1}^{n}(\alpha_j+\hat\alpha_j)\mu_j.
\end{eqnarray*}
\end{definition}

\begin{definition}
For the given weight $\mu=(\mu_1, \cdots, \mu_n)$, the \emph{weighted degree} of a monomial $A_{\alpha,\hat \alpha}z^\alpha \bar z^{\hat{\alpha}}$ is $l$, 
where $A_{\alpha,\hat \alpha} \in \mathbb C\setminus \{0\}$  and $z=(z_1, \cdots, z_n) \in \mathbb{C}^{n}$, if $|(\alpha,\hat \alpha)|_{\mu}=l$. 
A homogeneous polynomial $P$ is called \emph{$\mu$-homogeneous of weighted degree $l$},  or simply \emph{weighted homogeneous polynomial of degree $l$} , 
if it is a sum of monomials with weighted degree $l$.

 On the other hand, a vector field $Y$ is called \emph{$\mu$-homogeneous of weighted degree $l$}, or simply \emph{weighted vector field of degree $l$}, if it is a sum of 
 vector fields of the form
\begin{equation}
f(z,\bar{z})\partial_{z_j} \text{ or } f(z,\bar{z})\bar\partial_{z_j},
\end{equation}
where $f(z,\bar{z})$ is a weighted homogeneous polynomial of degree $l+\mu_j$.
Moreover, with respect to coordinates $(z, w)=(z_1, \cdots, z_n, w)$ of $\mathbb{C}^{n+1}$ with the weight $(\mu,1)=(\mu_1, \cdots, \mu_n,1)$,
a holomorphic vector field $Y$ is called \emph{$\mu$-homogeneous of weighted degree $l$}, or simply \emph{weighted holomorphic field of degree $l$}, 
if it is sum of holomorphic vector fields with polynomial coefficients, which are of the form
\begin{equation}
f(z,w)\partial_{z_j} \text{ and } g(z,w)\partial_{w},
\end{equation}
where $f(z,w)$ is a $(\mu,1)$-homogeneous of weighted degree $l+\mu_j$ and $g(z,w)$ is a $(\mu,1)$-homogeneous of weighted degree $l+1$.
\end{definition}

\begin{remark}\label{decompo} If $P$ is a $\mu$-homogeneous polynomial of weighted degree 1, we can write $P$ as
\begin{eqnarray}
P(z,\bar z) = \sum_{|(\alpha,\hat \alpha)|_{\mu}=1}A_{\alpha,\hat\alpha}z^{\alpha}\bar{z}^{\hat\alpha}.
\end{eqnarray}
\end{remark}

\begin{remark} The weighted degree of the vector fields $\partial_{z_j}$ and $\bar{\partial}_{z_j}$ are $-\mu_j$. The weighted degree of the vector field $\partial_{w}$ is $-1$.
\end{remark}

\bd Let $M$ be a hypersurface of $\mathbb{C}^{n+1}$ and $p \in M$ be a point. The weight $\mu$ is called \emph{distinguished} if there exist local holomorphic coordinates
$(z,w)$ centered at $p \in M$, where $z=(z_1, \cdots, z_n) \in \mathbb{C}^{n}$ and $w=u+iv$, such that $M$ is described as
\begin{eqnarray}\label{hyp}
v=P(z,\bar z)+\text{\it higher order terms}
\end{eqnarray}
at $p$, where $P(z,\bar z)$ is a weighted homogeneous polynomial of degree 1 without pluriharmonic terms. For the distingushed weight $\mu$, the local holomorphic 
coordinates are called $\mu$-adapted. 
Let $\mu_{M}=(\mu_1, \cdots, \mu_n)$ be the infimum of distingushed weights with respect to lexicographic ordering. The \emph{mulitype} of $M$ at $p$ is 
the $n$-tuple $(m_1,m_2,\cdots, m_n)$, where $m_j=\frac{1}{\mu_j}$ if $\mu_j\neq 0$ or $m_j=\infty$ if $\mu_j=0$. If $\mu_j\neq 0$ for all $j$, we say that $M$ is 
\emph{finite mulitype} at $p$. The weight $\mu_{M}$ is called the \emph{mulitype weight} and $\mu_{M}$-adapted coordinate are called \emph{mulitype coordinate}.
\ed

We will define a \emph{model hypersurface associated to $M$ at $p$} of the form \eqref{hyp}. $M_P$ is the hypersurface of $\mathbb{C}^{n+1}$ given by
\begin{eqnarray} \label{m-hyp}
\Im (w) = P(z,\bar z)
\end{eqnarray}
in the coordinates $(w,z)$, where $P(z,\bar z)$ is the weighted homogeneous polynomial of degree 1 with respect to the weight $\mu_M$, and we will
assume the hypersurface is pseudoconvex.

As an important example, let $P_j$ be the weighted homogeneous polynomials in $z=(z_1, \cdots, z_n) \in \mathbb{C}^{n}$ of weight $\frac{1}{2}\mu$,
where $\mu$ is a distinguished weight, and let $M_{P}$ be the hypersurface of $\mathbb{C}^{n+1}$ given by
\begin{eqnarray}
\Im (w) = \sum_{j=1}^{k} |P_j(z)|^2.
\end{eqnarray} If we assume that $deg\;  P_j \geq 1$, 
then $M_{P}$ is automatically pseudoconvex and we call it a \emph{sum of squares-type model}.

Let $w=u+iv$ and let $W$ be the vector field of degree $-1$ given by
\begin{eqnarray*} \partial_w =\partial_u -i \partial_v,
\end{eqnarray*}
then we have
\begin{eqnarray}
\partial_u( Im (w) - P(z,\bar z) ) =0
\end{eqnarray}
which gives us the \emph{first symmetry}.
As was proved in \cite{KM14}, the Lie algebra $\frak g$ of the infinitesimal automorphism $\aut(M,0)$ at $0 \in M \subset \mathbb{C}^{n+1}$ of $M$ admits a weighted decomposition:
\begin{eqnarray}
\frak g =\frak g_{-1} \oplus \bigoplus_{j=1}^{n} \frak g_{-\mu_j} \oplus \frak g_0 \oplus \frak g_c \oplus \frak g_n \oplus \frak g_1.
\end{eqnarray}
where the local vector field in $\frak g_c$ commute with $W$, the non-zero local vector field in $\frak g_n$ do not commute with $W$ and their weights $\mu_j$ are between 0 and 1. In particular, $W = \partial_w$ is contained in $\frak g_{-1}$, 
which has real dimension one.

Let $E$ be the \emph{Euler field} given by
\begin{eqnarray}
w \partial_w + \sum_{j=1}^{n}\mu_j z_j \partial_{z_j}
\end{eqnarray}
then it is immediate that
\begin{eqnarray}
\Re E( Im (w) - P(z,\bar z) ) = 0
\end{eqnarray}
which implies $E \in \frak g_0$ and hence, $dim(\frak g_0)\geq 1$. Thus the Euler field $E$ gives the \emph{second symmetry} in $\aut (M_P, 0) $, for an arbitrary model $M_P$.

\section{Rotations}

Let $X$ be a weight zero infinitesimal CR automorphism in $\mathfrak g_0$. By the result of \cite{KM14}, $X$ is a linear vector field. It's Jordan normal form could be decomposed into $X^{Re}+X^{Im}+X^{Nil}$, where $X^{Re}$ is the real diagonal of the Jordan normal form, $X^{Im}$  is the imaginary diagonal of the Jordan normal form and  $X^{Nil}$ is nilpotent part. In particular, if we let $\tilde E$ be the vector field given by
\begin{eqnarray}
i\sum_{j=1}^{n}\mu_j z_j \partial_{z_j},
\end{eqnarray} then $\tilde E =\tilde E^{Im}$.

\bl \label{diagonal} Let the model hypersurface \eqref{m-hyp} be  pseudoconvex and consider the expanded polynomial of Remark \ref{decompo}.
There exist multitype coordinates in which for each $j = 1, \dots , n$ 
there exists a multiindex $\alpha$ such that $A_{\alpha, \alpha} \neq 0$, 
$\alpha_j \neq 0$ and  $\alpha_l = 0$ for all $l>j$. 
 \el

\begin{proof}
By Lemma 3.1 of \cite{K10}, there exist multitype coordinates such that for each $j$ 
the j-th partial derivative of the polynomial 
$$ P^j(z_1, \dots, z_j) := P (z_1, \dots, z_j, 0,0,  \dots, 0)$$
is not identically zero, i.e. $$\frac{\partial P^j}{\partial z_j} \neq 0.$$
For any choice of the variables $z_1, \dots, z_{j-1}$, consider the 
restriction of $P_j$ to the complex line $z_1=c_1, \dots, z_{j-1}=c_{j-1}$ and denote this one variable polynomial by $Q$. 
Since $\Delta Q$ is nonnegative, it must contain a leading term of the form $\vert z_j \vert^{2m}$, of the lowest degree. We will consider the minimal
such $m$ over the choice of the complex lines $z_1=c_1, \dots, z_{j-1}=c_{j-1}$. Then $P^j$ contains the summand $S(z_1, \dots, z_{j-1})\vert z_j \vert^{2m}$. $S$ has to be a 
nonnegative polynomial, hence it contains a nonzero term with $\alpha = \hat \alpha$  in $z_1, \dots, z_{j-1}$. That finishes the proof.
\end{proof}

Let us remark that by Lemma 4.6. in \cite{KM14}, if $X$ lies in $\aut(M_P, 0)$, then in Jordan normal form both its diagonal and nilpotent part 
lie in $\aut(M_P, 0)$. Moreover, both the real and imaginary parts of the diagonal component also lie in   $\aut(M_P, 0)$.

\bl \label{real} Let $X$ be a weight zero rigid infinitesimal CR automorphism in $\mathfrak g_0$ and $X^{Re}$ be the real diagonal of the Jordan normal form of $X$. 
Then $X^{Re}=0$ 
\el

\begin{proof}
 Let $X =X^{Re}$ be the vector field
\begin{eqnarray}
\sum_{j=1}^{n}\lambda_j z_j \partial_{z_j},
\end{eqnarray} where $\lambda_j \in \mathbb R$, and let $z^{\alpha}\bar{z}^{\hat\alpha}$ be the monomial, where $z=(z_1, \cdots, z_n)$, $\alpha=(\alpha_1,\cdots \alpha_n)$ and $\hat\alpha=(\hat\alpha_1,\cdots \hat\alpha_n )$. We apply $X$ and $\bar X$ to $z^{\alpha}\bar{z}^{\hat\alpha}$, then
\begin{eqnarray*}
X(z^{\alpha}\bar{z}^{\hat\alpha})=\sum_{j=1}^{n}\lambda_j \alpha_j z^{\alpha}\bar{z}^{\hat\alpha}, \hspace{.5cm}
\bar{X}(z^{\alpha}\bar{z}^{\hat\alpha})=\sum_{j=1}^{n}\lambda_j \hat\alpha_j z^{\alpha}\bar{z}^{\hat\alpha}
\end{eqnarray*}
and hence,
\begin{eqnarray*}
0=2ReX(z^{\alpha}\bar{z}^{\hat\alpha})=\sum_{j=1}^{n}\lambda_j (\alpha_j+\hat\alpha_j) z^{\alpha}\bar{z}^{\hat\alpha}
\end{eqnarray*}
if and only if
\begin{eqnarray}
\sum_{j=1}^{n}\lambda_j (\alpha_j+\hat\alpha_j)=0.
\end{eqnarray}
We denote this $\lambda \perp (\alpha+\hat\alpha)$.

Let $P$ be the homogeneous polynomial, then we have decomposition based on the monomials
\begin{eqnarray}
P(z,
\bar z) = \sum_{\alpha, \hat\alpha}A_{\alpha,\hat\alpha}z^{\alpha}\bar{z}^{\hat\alpha}.
\end{eqnarray}
Then $X P(z,
\bar z)$ and $\bar{X}P(z,
\bar z)$ are
\begin{eqnarray*}
XP(z,
\bar z) 
&=& \sum_{\alpha, \hat\alpha}\sum_{j=1}^{n}\lambda_j \alpha_j A_{\alpha,\hat\alpha}z^{\alpha}\bar{z}^{\hat\alpha} \\
\bar{X}P(z,
\bar z) 
 &=& \sum_{\alpha, \hat\alpha}\sum_{j=1}^{n}\lambda_j \hat\alpha_j A_{\alpha,\hat\alpha}z^{\alpha}\bar{z}^{\hat\alpha}.
\end{eqnarray*}
Hence,
\begin{eqnarray*}
0=2ReXP(z,
\bar z)=\sum_{\alpha, \hat\alpha}\sum_{j=1}^{n}\lambda_j (\alpha_j+\hat\alpha_j)A_{\alpha,\hat\alpha}z^{\alpha}\bar{z}^{\hat\alpha}
\end{eqnarray*}
if and only if $\lambda \perp (\alpha+\hat\alpha)$ for all $A_{\alpha,\hat\alpha} \neq 0$.

By Lemma \ref{diagonal}, the subset $\{2\alpha \in \mathbb C^n| A_{\alpha,\alpha} \neq 0\}$ of $\{\alpha+\hat\alpha \in \mathbb C^n| A_{\alpha,\hat\alpha} \neq 0\}$ spans $\mathbb C^n$. Hence, $\lambda=0$.
\end{proof}


\bl \label{nil} Let $X$ be a weight zero rigid infinitesimal CR automorphism in $\mathfrak g_0$ and $X^{Nil}$ be the nilpotent part of the Jordan normal form of $X$.
Then $X^{Nil}=0$ \el

\begin{proof}
 Let $Y=X^{Nil}$ be the nilpotent part of $X$, 
\begin{eqnarray}
\sum_{i=1}^{n-1}\lambda_i z_i \partial_{z_{i+1},}
\end{eqnarray} where $\lambda_i \in \{0,1\}$.

We expand $P$ based on the monomials;
\begin{eqnarray*}
P(z,\bar{z})&=& \sum_{\alpha, \hat\alpha}A_{\alpha,\hat\alpha}z^{\alpha}\bar{z}^{\hat\alpha}.
\end{eqnarray*}
From the equation $Y P(z,\bar{z})= - \overline{Y} P(z,\bar{z})$, we have
\begin{eqnarray}
\label{conj-symmetry} \sum_{\alpha, \hat\alpha}\sum_{i=1}^{n}\lambda_i\alpha_{i+1}A_{\alpha,\hat\alpha}z^{\alpha(i)}\bar{z}^{\hat\alpha}&=&-
\sum_{\alpha, \hat\alpha}\sum_{i=1}^{n} \lambda_i\hat\alpha_{i+1}A_{\alpha,\hat\alpha}z^{\alpha}\bar{z}^{\hat\alpha(i)},
\end{eqnarray}
 where $\alpha(i)=(\alpha_1, \cdots, \alpha_i+1, \alpha_{i+1}-1, \cdots, \alpha_n)$ and $\hat\alpha(i)=(\hat\alpha_1, \cdots, \hat\alpha_i+1, \hat\alpha_{i+1}-1, \cdots, \hat\alpha_n)$.

For the set $\{\alpha\}=\{\hat\alpha\}$, we give a partial ordering by $\alpha \succ_1 \beta$ if $\alpha(i)=\beta$ for some $i \in \{1,\cdots, n-1\}$: we denote $\alpha \succeq \beta$ if $\alpha = \beta$, $\alpha \succ_1 \beta$ or there is a subset $\{\gamma_i|i=1,\cdots t, \gamma_1\succ_1, \cdots, \succ_1\gamma_t \}$ of $\{\alpha\}$ such that $\alpha \succ_1 \gamma_1$ and $\gamma_t  \succ_1 \beta$. Then $\{\alpha, \succeq \}$ is a partially ordered set.

 Assume $Y \neq 0$. Let $i$ be the largest integer with $\lambda_i \neq 0$. By Lemma \ref{diagonal}, there exists $\alpha$ such that $\alpha_{i+1}\neq 0$ and $A_{\alpha,\alpha} \neq 0$. In the equation (\ref{conj-symmetry}), non-vanishing $\lambda_i\alpha_{i+1}A_{\alpha,\alpha}$ of the left-side implies there are $\beta$ and $\hat\beta$ such that $\lambda_k\hat\beta_{k+1}A_{\beta,\hat\beta}\neq 0$ for some $k$ and \begin{eqnarray*} z^{\alpha(i)}\bar{z}^{\alpha} = z^{\beta}\bar{z}^{\hat\beta(k)} \end{eqnarray*}
which gives us $\hat\beta \succ_1 \alpha$.
Since $\lambda_k\hat\beta_{k+1}\neq 0$, the equation (\ref{conj-symmetry}) implies there are $\gamma$ and $\hat\gamma$ such that $\lambda_k\hat\gamma_{l+1}A_{\beta,\hat\beta}\neq 0$ for some $l$ and \begin{eqnarray*} z^{\hat\beta(k)}\bar{z}^{\hat\beta}= z^{\gamma}\bar{z}^{\hat\gamma(l)} \end{eqnarray*}
implies $\hat\gamma \succ_1 \hat\beta \succ_1 \alpha$. If we continue this, we get infinite-length chain, which is a contradiction, since the set $\{\alpha\}$ is finite. Hence, $Y=0$.
\end{proof}

\section{Generalized rotations}

In this section we prove the following result.

\bt
\label{gc}
If $M_P$ given by \eqref{m-hyp} is pseudoconvex, then $\frak g_c = 0$. 
\et

\begin{proof}
Since $M_P$ is pseudoconvex near $0$, any neighbourhood of $0$ contains strongly pseudoconvex points, 
by the finite multitype assumption. It follows that $\aut(M_P,0)$ is a subalgebra of $\aut(Q,0)$, where $Q$ is the strongly pseudoconvex hyperquadric.
We know that  $\aut(Q,0)$ admits the weighted grading 

\begin{eqnarray}
\frak g =\frak g_{-1} \oplus \frak g_{-\frac12} \oplus \frak g_0 \oplus \frak g_{\frac12} 
\oplus \frak g_1
\end{eqnarray}
Here both $ \frak g_{-1}$ and $\frak g_{1}$ are one dimensional. Moreover, $\frak g_0$ can be split further as $\frak g_0 = \frak g_{0e} \oplus \frak g_{0r}$, where
$\frak g_{0e} $ is the one dimensional subspace contaning the Euler field and  $\frak g_{0r}$
contains rigid vector fields, i.e., rotations. 
We have the following commutation relations  $$[\frak g_{-1}, \frak g_{1}] = \frak g_{0e}, \ \ \ \ \ \ \  
[\frak g_{0e}, \frak g_{0r}] = 0.$$
Moreover,  $\frak g_{0r} $ commutes both with  $\frak g_{-1}$ and $\frak g_{1}$. 

Note that $\frak g_{0e}$ is spanned by the Euler field $E=E_Q$, which is the grading element. Let us emphasize the crucial fact that on the strongly psedoconvex hyperquadric 
the Euler field is determined uniquely, since there exist no real rotations.

Assume  $\dim \frak g_c > 0$ and let  $Z \in \frak g_c$. Denote by   $\frak h$ the subalgebra of 
 $\frak g$ generated by the vector fields $E, W, Z$. $\frak h$ is three dimensional with the commutation relations
 \begin{equation}\label{ew1}
[E, W] = -W,
\end{equation}
 \begin{equation}\label{ez1}
[E, Z] = cZ
\end{equation}
where $ 0 < c < 1$, and 
 \begin{equation}\label{zw1}
[Z, W] = 0.
\end{equation}
We will prove that $\aut(Q, 0)$ contains no subalgebra isomorphic to $\frak h$. 
By contadiction, asuume that  $\aut(Q, 0)$ contains such an algebra and denote the corresponding vector fileds in  $\aut(Q, 0)$ 
by primes. Hence  $E', W', Z'$ generate a three dimensional subalgebra with the same relations as in \eqref{ew1} - \eqref{zw1}.

We will decompose  $E', W', Z'$ into the graded components of  $\aut(Q, 0)$ and then express the commutators of those 
three vector fileds according to this grading. More precisely, we will be only interested in the $\frak g_{-1}$ $\frak g_{0e}$ 
and 
$\frak g_{1}$ components of these commutators.

Note that for this purpose  we may ignore the $\frak g_{-\frac12}$ and the $\frak g_{\frac12}$ components of $\frak g$. 
Indeed,  $[\frak g_{-\frac12}, \frak g_{\frac12}] = 0$. and the nonvanishing commutators of  $\frak g_{-\frac12}$ and  $\frak g_{\frac12}$ 
with the components $\frak g_{-1}$, $\frak g_{0}$ and $\frak g_{1}$  lie in  $\frak g_{-\frac12}$ or $\frak g_{\frac12}$.

Hence we can write

\begin{equation}\label{zw}
E' = e_1 W_Q + e_2 E_Q + e_3 G_1 + \epsilon_E
\end{equation}

\begin{equation}\label{zw}
W' = w_1 W_Q + w_2 E_Q + w_3 G_1 + \epsilon_W
\end{equation}

\begin{equation}\label{zw}
Z' = z_1 W_Q + z_2 E_Q + z_3 G_1 + \epsilon_Z
\end{equation}
where  $W_Q \in \frak g_{-1}$ $E_Q \in \frak g_{0e}$ and $G_1 \in \frak g_{1}$ satisfying equation (\ref{ew})-(\ref{zw}) and $\epsilon_E$, $\epsilon_W, \epsilon_Z$, represent terms which will not contribute to the commutator 
components on the level of $\frak g_{-1}$, $\frak g_{0e}$ and $\frak g_{1}$, in other words they 
 belong to $\frak g_{-\frac12} \oplus \frak g_{0r} \oplus \frak g_{\frac12}$. 
 
We have
 \begin{equation}\label{ew}
[E_Q, W_Q] = -W_Q,
\end{equation}
 \begin{equation}\label{ez}
[E_Q, G_1] = G_1
\end{equation}
and 
 \begin{equation}\label{zw}
[G_1, W_Q] = -E_Q.
\end{equation}

Hence we obtain  for $[E',W']$

\begin{eqnarray*}
e_2 w_1 - e_1 w_2 &=  w_1 \\
e_1 w_3 - e_3 w_1 &=  w_2 \\
e_2 w_3 - e_3 w_2 &=  w_3. 
\end{eqnarray*}

For $[E',Z']$ we have 
 
\begin{eqnarray*}
e_2 z_1 - e_1 z_2 = c z_1 \\
e_1 z_3 - e_3 z_1 = c z_2 \\
e_2 z_3 - e_3 z_2 = c z_3, 
\end{eqnarray*}
and for 
$[Z',W']$ we obtain
 
\begin{equation} \label{3rd}
z_2 w_1 - z_1 w_2 = 
z_1 w_3 - z_3 w_1 =
z_2 w_3 - z_3 w_2 = 0. 
\end{equation}

By the third set of equations \eqref{3rd}, the vector product of $z$ and $w$ is zero, hence $z$ and $w$ are linearly dependent. 
From the first two sets we obtain immediate contradiction.

\end{proof}

\section{Sum of squares models}

In this section we consider sum of squares models and give a description of their symmetry algebra.

\bt
Let  $M_{S}$ be the sum of squares  homogeneous polynomial model of degree $ k > 2$. Then the subspace $\frak g_{-\mu_j}$ $j=1,..n$ of the symmetry algebra vanishes.
\et

\begin{proof} Since the local vector field $\partial_{z_i}$ has the weight $-\mu_i$, it is equivalent to show that \begin{eqnarray}
Re\partial_{z_i}( Im (w) - \sum_{j=1}^{k} |P_j(z)|^2 ) \neq 0.
\end{eqnarray}
Since $\partial_{z_i}Im (w)=0$, it is enough to prove $Re\partial_{z_i} \sum_{j=1}^{k} |P_j(z)|^2 \neq 0$. Since $M_S$ is holomorphically non-degenerate,
i.e., $\{\bigtriangledown P_j\}$ spans $\mathbb C^n$ at a generic point, we have $max_{j}\{deg_{z_i} P_j(z)\} > max_{j}\{deg_{z_i} \partial_{z_i}P_j(z)\}\geq 0$. It follows that
\begin{eqnarray}
\sum_{j=1}^{k} (\partial_{z_i}P_j)\overline{P_j} \neq -\sum_{j=1}^{k}(\overline{\partial_{z_i}P_j}){P_j}.
\end{eqnarray}
\end{proof}

Let $X =X^{Im}$ be the vector field
\begin{eqnarray}
\sum_{j=1}^{n}i\lambda_j z_j \partial_{z_j},
\end{eqnarray} where $\lambda_j \in \mathbb R$, and let $z^{\alpha}\bar{z}^{\hat\alpha}$ be the monomial,
where $z=(z_1, \cdots, z_n)$, $\alpha=(\alpha_1,\cdots \alpha_n)$ and $\hat\alpha=(\hat\alpha_1,\cdots \hat\alpha_n )$. 
We apply $X$ and $\bar X$ to $z^{\alpha}\bar{z}^{\hat\alpha}$, then
\begin{eqnarray*}
X(z^{\alpha}\bar{z}^{\hat\alpha})=\sum_{j=1}^{n}i\lambda_j \alpha_j z^{\alpha}\bar{z}^{\hat\alpha}, \hspace{.5cm}
\bar{X}(z^{\alpha}\bar{z}^{\hat\alpha})=\sum_{j=1}^{n}-i\lambda_j \hat\alpha_j z^{\alpha}\bar{z}^{\hat\alpha}
\end{eqnarray*}
and hence,
\begin{eqnarray*}
0=2ReXP(z,\bar{z})=\sum_{\alpha, \hat\alpha}\sum_{j=1}^{n}i\lambda_j (\alpha_j-\hat\alpha_j)A_{\alpha,\hat\alpha} z^{\alpha}\bar{z}^{\hat\alpha}
\end{eqnarray*}
if and only if
\begin{eqnarray}
\sum_{j=1}^{n}\lambda_j\alpha_j=\sum_{j=1}^{n}\lambda_j\hat\alpha_j.
\end{eqnarray}

\bl \label{e tilde} If $P(z,\bar z)= \sum_{j=1}^{k} |P_j(z)|^2$ is a sum of squares polynomial of weighted degree 1, then $\tilde E \in \frak g_0$.
\el

\begin{proof}
If $P(z,\bar z)= \sum_{j=1}^{k} |P_j(z)|^2$ is a sum of squares-type, then
\begin{eqnarray}
\sum_{j=1}^{n}\mu_j\alpha_j=\sum_{j=1}^{n}\mu_j\hat\alpha_j=\frac{1}{2}
\end{eqnarray}
which implies the result.
\end{proof}

In this case: if $P(z,\bar z)= \sum_{j=1}^{k} |P_j(z)|^2$ is a sum of square-type of degree $l$, by Theorem 4.7 of \cite{KM14} and by Lemma \ref{e tilde}, the vector field
\begin{eqnarray}
\frac{1}{2}w^2 \partial_w + 2\sum_{j=1}^{n}\mu_j w z_j \partial_{z_j}
\end{eqnarray}
is contained in $\mathfrak g_1$. It provides the \emph{third symmetry}.

As a consequence, we obtain the following precise description of the Lie algebra of infinitesimal automorphisms of $M_S$.

\bt
Let  $M_{S}$ be the sum of square-type homogeneous polynomial model $M_{S}$ of degree $ k > 2$. 
Then the   Lie algebra of infinitesimal automorphisms
$\5g=\aut(M_{S},0)$ of $M_{S}$ admits the weighted grading given
by
\begin{equation}
\5g = \5g_{-1} \oplus \5g_{ 0} \oplus \5g_{1}
\end{equation}
where $\5g_{-1}$ and $\5g_{1}$ are of real dimension one and $\5g_{0}$ is generated by the Euler field and a subalgebrs of $\mathfrak u(n)$. 

\et

Now we consider the general case of a weighted homogeneous polynomial model. We will detone by $\kappa_M$ the number of multitype weights with 
$\mu_j = \frac12$. 

\bt
Let  $M_{S}$ be the sum of square-type weighted homogeneous polynomial model $M_{S}$ of degree $ k > 2$. 
Then the   Lie algebra of infinitesimal automorphisms
$\5g=\aut(M_{S},0)$ of $M_{S}$ admits the weighted grading given
by
\begin{equation}
\5g = \5g_{-1} \oplus \5g_{-\frac12} \oplus\5g_{ 0} \oplus\5g_{\frac12 } \oplus \5g_{1}
\end{equation}
where $\5g_{-1}$ and $\5g_{1}$ are of real dimension one,  $\5g_{-\frac12}$ and  $\5g_{-\frac12} $ are of real dimension $2 \kappa_M$,  and $\5g_{0}$ is generated by the Euler field and a subalgebra of $\mathfrak u(n)$. 

\et

\begin{proof}
 It follows from the assumptions that in suitable multitype coordinates we can write 
 \begin{equation}
  P(z, \bar z) = \sum_{j=1}^{ \kappa_M} \vert z_j \vert^2 +  Q(z_{\kappa_M +1}, \dots, z_n, \bar z_{\kappa_M +1}, \dots, \bar z_n)
\end{equation}
where $Q$ is weighted homogeneous and balanced.
We verify that the vector fields 

\begin{eqnarray}
a  \partial_{z_j} + 2i \bar a z_j \partial_{w}
\end{eqnarray}
for each $j =1, \dots, \kappa_M$ and $a \in \mathbb C$ lie in  $\5g_{-\frac12}$, and
can be integrated to  vector fields 
\begin{eqnarray}
a w \partial_{z_j} + 2i \bar a z_j \sum_{k=1}^n \mu_k z_k \partial_{z_k}
+ 2i \bar a z_j w \partial_{w}
\end{eqnarray}
which lie in  $\5g_{\frac12}$. By the above Theorem 5.1,
there are no other elements of  $\5g_{-\frac12}$, hence no other elements in  $\5g_{\frac12}$, 
which gives the claim. 
\end{proof}

\section{Proof of the main results}
We give now the proof of Theorem 1.1.

\begin{proof}
 By Theorem 1.3 of \cite{KM14} and Theorem 4.1, we obtain that 
 $\5g=\aut(M_{P},0)$ admits the weighted grading given
by

\begin{eqnarray}
\frak g =\frak g_{-1} \oplus \bigoplus_{j=1}^{n} \frak g_{-\mu_j} \oplus \frak g_0 \oplus \frak g_n \oplus \frak g_1.
\end{eqnarray}

It remains to prove that the weight of $g_n$ is $\frac12$. We will use the following characterization of manifolds with nonvanishing $\5g_n$. 

Let $M_P$ have nontrivial $\5g_n$
and
\begin{equation}
  X=i \partial _ {z_{l}}
  \label{xi}
\end{equation}
for some $l$ be an infinitesimal symmetry of $M_P$, which can be integrated.
Let us write $P$ as 
\begin{equation}
{P}(z, \bar z) = \sum_{j=0}^{m}{(\Re z_l)}^j P_j(z', \bar z'),\ \
 \label{pzz}
\end{equation}
for some homogeneous 
polynomials  ${P_j}$
in the variables $z' = ( z_1, \dots, \hat z_{l},  \dots, z_n)$,
with $P_m \neq 0.$
Then $M_P$  has one of the following two forms. 
Either 
\begin{equation}
 P(z, \bar z) =
 x_l^2 + x_l P_1( z', \bar z') +
 P_0( z', \bar z'),
\end{equation}
or 
\begin{equation}
 P(z, \bar z) =
 x_l P_1( z', \bar z') +
 P_0( z', \bar z').
\end{equation}


This characterization was proved 
by Kol\'a\v r and Meylan in \cite{KM18} in the $\mathbb C^3$ case, and in \cite{KM19} in general.
It is immediate to verify that in the second case, the manifold is not pseudoconvex. 
On the other hand, in the first case we have $\mu_l = \frac12$, which leads to the claim of Theorem 1.1.

 \end{proof}

Theorem 1.2 is an immediate consequence of Theorem 1.1

\end{document}